\documentclass[12pt]{article}
\pdfoutput=1 
\usepackage{graphicx}
\usepackage{bm,amsmath,amssymb,mathrsfs}
\usepackage{amsfonts}
\usepackage{algorithm,algpseudocode,float}
\usepackage{algorithmicx}
\usepackage{multirow}

\makeatletter
\newenvironment{breakablealgorithm}
  {
   \begin{center}
     \refstepcounter{algorithm}
     \hrule height.8pt depth0pt \kern2pt
     \renewcommand{\caption}[2][\relax]{
       {\raggedright\textbf{\ALG@name~\thealgorithm} ##2\par}%
       \ifx\relax##1\relax 
         \addcontentsline{loa}{algorithm}{\protect\numberline{\thealgorithm}##2}%
       \else 
         \addcontentsline{loa}{algorithm}{\protect\numberline{\thealgorithm}##1}%
       \fi
       \kern2pt\hrule\kern2pt
     }
  }{
     \kern2pt\hrule\relax
   \end{center}
  }
\makeatother

\makeatletter
\let\OldStatex\Statex
\renewcommand{\Statex}[1][3]{%
  \setlength\@tempdima{\algorithmicindent}%
  \OldStatex\hskip\dimexpr#1\@tempdima\relax}
\makeatother

\def\p{\partial}
\def\DD{\displaystyle}

\def\F{{\mathcal{F}}}
\def\G{{\mathcal{G}}}
\def\U{{\tilde{U}}}
\def\l{\,l}

\def\a{\alpha}
\def\lpml{l_{\text{pml}}}
\def\FF{\text{F}}
\def\II{\text{I}}
\def\lx{\bar{l}_1}

\def\L{{\mathscr{L}}}

\def\FSUP{}

\title{A Fast Propagation Method \\
  for the Helmholtz equation}

\author{Wei Leng \footnote{State Key Laboratory of Scientific and Engineering Computing,
Chinese Academy of Sciences, Beijing 100190, China. Email: {wleng@lsec.cc.ac.cn}.}}

\begin{document}

\maketitle

\abstract{A fast method is proposed for solving the high frequency 
  Helmholtz equation.
   The building block of the new fast method is an overlapping source transfer domain decomposition method for layered medium, which is an extension of the source transfer domain decomposition method proposed by Chen and Xiang \cite{Chen2013a,Chen2013b}.
  The new fast method contains a setup phase and a solving phase. In the
  setup phase, the computation domain is decomposed hierarchically into 
  many subdomains of different levels,
  and the mapping from incident traces to field traces on all the subdomains
  are set up bottom-up.   
  In the solving phase, first on the bottom level, the local problem
  on the subdomains with restricted source is solved, then the wave 
  propagates on the boundaries of
  all the subdomains bottom-up, at last the local solutions on all the subdomains
  are summed up top-down.
  The total computation cost of the new fast method is $O(n^{\frac{3}{2}} \log n)$ for 2D problem.
  Numerical experiments shows that with the new fast method, Helmholtz
  equations with half billion
  unknowns could be solved efficiently on massively parallel
  machines.
}



\section{Introduction}

We consider in this paper to solve the Helmholtz equation in the full space $\mathbb{R}^2$, with Sommerfeld radiation condition,
\begin{align} \label{eq:helm}
  \Delta u + k^2 u &= f \qquad  \mbox{in} \,\,\, \mathbb{R}^2, \\
  r^{1/2} (\frac{\p u}{\p r} - \mathbf{i} k u) &\rightarrow 0 \qquad \mbox{as} \,\,\, r = |x| \rightarrow \infty \nonumber 
\end{align}
where $k$ is the wave number.

%
Many domain decomposition method has recently been developed to solve the Helmholtz equation, 
most of them are non-overlapped, and the major differences are the interface
conditions.
Engquist and Ying \cite{Engquist2011a, Engquist2011b} proposed a sweeping preconditioner
by approximating the inverse of Schur complements in the LDL$^t$ factorization,
Stolk \cite{Stolk2013} proposed a domain decomposition method with a transmission condition based on perfect matched layers, 
Vion an Geuzaine \cite{Vion2014} proposed a double sweep preconditioner that
use a transmission condition that involves Dirichlet-to-Neumann (DtN) operator, 
Zepeda \cite{Zepeda2014} introduced the method of polarized trace that use a
transmission condition in boundary integral form,  
Liu and Ying \cite{Liu2015} developed an additive sweeping preconditioner that 
use a transmission condition built with the boundary values of the intermediate 
wave directly. 
Chen and Xiang \cite{Chen2013a,Chen2013b} proposed the source transfer domain
decomposition method that transfer the source in subdomains,
and recently Du and Wu \cite{Wu2015} improved the method so that the transfer applies in both directions.

The domain decomposition method in the literature usually 
approximately solves the Helmholtz equation with varying medium, either with approximated interface
condition or with approximated Green function, thus they are commonly used as preconditioners for Krylov subspace method such as GMRES.

An overlapping source transfer domain decomposition method is proposed 
for Helmholtz equation with layered medium, the method follows the natural 
wave traveling process in layered medium, which involves the
reflections and refractions at the interface of the layers. The convergence 
of the new domain decomposition method is ensured by the overlapping region, 
and the accuracy of the new domain decomposition method makes it
 the building block of the new fast method.

The domain decomposition method suffers from slow convergence rate when the 
number of subdomains is large,
thus multilevel grid is needed so that the information is brought to far away subdomains
  without passing the subdomains on the way.
The upper level grid for Poisson type problem could be coarser since 
the amount of information decreases fast as the distance grows. 
However, 
for Helmholtz equation, the grid size should be maintained small
to represent wave shapes on the upper level grid.
Fortunately, the trace on the subdomain boundaries could be used to 
represent the solution on the subdomain, thus the computation cost 
on upper level grid is not formidable. 

The fast method we proposed first setup the trace mapping on subdomains
of different levels. Then the sources are converted to traces on the bottom
level, and propagate on higher and higher levels till 
the top level,
 then the traces on high levels are decomposed into traces on
lower and lower level, at last the traces in 
the bottom level is converted back to solutions and summed up. 
In such up and down process, the wave travels to 
far away regions via the traces on high levels.

The rest of the paper is organized as follows. 
In section 2, an overlapping source transfer domain decomposition method is 
proposed for Helmholtz equation with layered medium.
In section 3, the fast algorithm is described. The multilevel domain decomposition
with quadtree structure is built, and the algorithm to build incident trace to field trace mapping on subdomains is proposed, 
then source up and solution down algorithm are proposed.
The numerical experiment for Marmousi model is present in section 4.

%
%
%
%

\section{The overlapping source transfer DDM}

The foundation of the fast method is the overlapping source transfer domain
decomposition method for the Helmholtz equation. 
We first propose and analyze the overlapping STDDM 
for Helmholtz problem with three layered medium, 
then revise the method and substitute the solving of subdomain problem into mapping,
and at last propose the overlapping decomposition method for four subdomains, 
which is the building block of the fast propagation method.

%
%
\subsection{STDDM in three layered medium}

\begin{figure}[ht!]
\centerline{
\includegraphics[width=.6\textwidth]{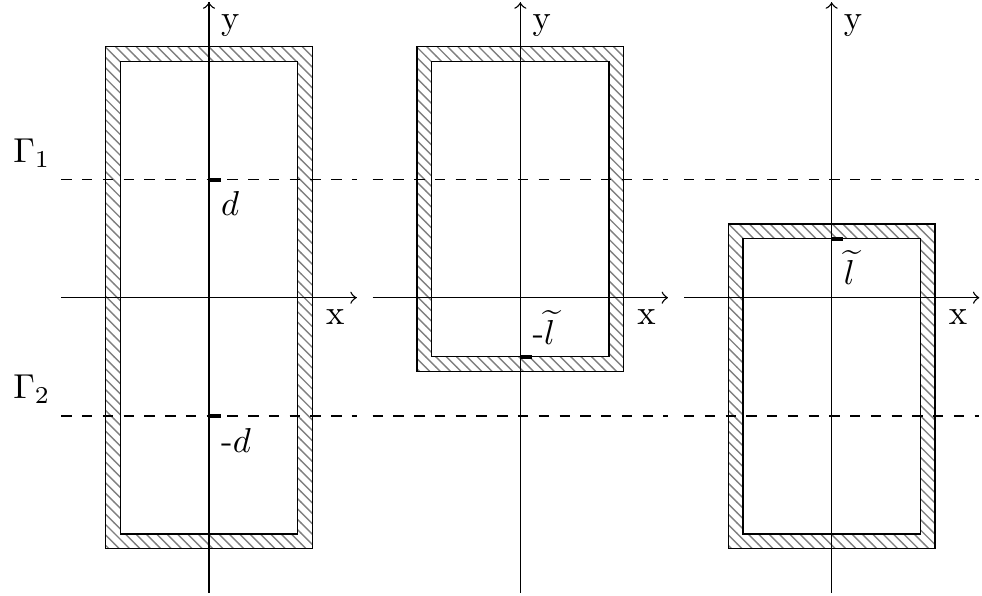}
} 
\caption{Domain decomposition in $y$ direction for three layered problem.} 
\label{fig:3layer}
\end{figure}

Consider the Helmholtz equation \eqref{eq:helm} defined in $\mathbb{R}^2$,
where the source $f$ is given, and the wave number $k$ is different in three horizontal layers,  
\begin{equation}
k(y) = \left\{
\begin{array}{ll}
 k_1,& \quad \mbox{if ~ $y < -d$}\\
 k_2,& \quad \mbox{if ~ $-d \le y \le d$}\\
 k_3,& \quad \mbox{if ~ $y > d$}\\
\end{array}
\right.
\end{equation}
as shown in Fig \ref{fig:3layer}. 
The upper interface $y=d$ 
is denoted $\Gamma_1$, and 
the lower interface $y=-d$ 
is denoted $\Gamma_2$.

The frequence domain wave equations defined on unbounded domain could
be solved on truncated domain with the perfect matched layer as the absorbing 
boundary condition \cite{Berenger, Chew1994}. 
To solve Helmholtz problem \eqref{eq:helm}, the unbounded domain $\mathbb{R}^2$ 
is truncated to a rectangle $[-l_1, l_1] \times [-l_2, l_2]$, 
with a PML layer of length $\lpml$ attached to the boundary, 
and the trucated domain $\Omega$
becomes  $[-l_1-\lpml, l_1+\lpml] \times [-l_2-\lpml, l_2+\lpml]$. 
We refer the domain without PML layer as the interior of domain $\Omega$, 
denoted $\tilde{\Omega}$. For simplicity , we denote $-l_1-\lpml$ as $\lx$.

The uniaxial PML method \cite{Chew1994} is used in this paper, 
where the complex coordinate is streched 
in $x$ and $y$ direction sperately,
$\DD \tilde{x}_j(x_j) = \int_0^{x_j} \sigma_j(t) dt$, $j = 1,2$, and
the medium perporty is chosen that $\sigma_j(t) = 0$ for $ |t| \leq l_j$, and $\sigma_j(t) > 0$
in PML layer $|t| > l_j$. Then the PML equation is
\begin{equation} \label{eq:pml}
  J^{-1} \nabla \cdot (A \nabla u) + k^2 u = f, \qquad \mbox{in} \,\,\, \Omega,
\end{equation}
where $\DD A(x) = \mbox{diag}\left(\frac{\a_2(x_2)}{\a_1(x_1)}, \frac{\a_1(x_1)}{\a_2(x_2)}\right)$, and $J(x) = \a_1(x_1) \a_2(x_2)$.

The computation domain is decomposed to two overlapping subdomains,
the upper one  $\Omega_1 = [-\lx, \lx] \times
[-\tilde{l}-\lpml, l_2+\lpml]$
and the lower one  $\Omega_2 = [-\lx, \lx] \times
[-l_2-\lpml, \tilde{l}+\lpml]$, with an overlapping region 
$[-\lx, \lx] \times [-\tilde{l}, \tilde{l}]$,
as is shown in Fig \ref{fig:3layer}.
Similar PML equations as \eqref{eq:pml} are built on the two subdomains,
and the parameter $A$ and $J$ in the PML equation are denoted $A_i$ and $J_i$
for subdomain $\Omega_i$, $i = 1,2$.

The new domain decomposition method first solve the subdomain problem 
with the restricted source,   
\begin{equation} \label{eq:dm1}
  J_i^{-1} \nabla \cdot (A_i \nabla u_i) + k^2 u_i = f_i, \qquad \mbox{in} \,\,\, \Omega_i, i = 1,2
\end{equation}
where $f_1 = f \cdot \chi_{y < 0}$ for $\Omega_1$, 
and $f_2 = f \cdot \chi_{y \geq 0}$ for $\Omega_2$, and the solution is denoted $u_i^0$ for $i = 1,2$.

Then, 
the wave field in $\Omega_1$ is transfered as source to $\Omega_2$ 
meanwhile the wave field in $\Omega_2$ is transfered as source to $\Omega_1$, 
 with the new transfered sources the PML equation on the subdomains
is solved and new wave field is generated,
and so on,
\begin{align} 
J_1^{-1} \nabla \cdot (A_1 \nabla u_1^{s+1}) + k^2 u_1^{s+1} = \Psi_1(u_2^s), & \qquad \mbox{in} \,\,\, \Omega_{1} \label{eq:tran1} \\
\Psi_1(u_2^s) = - J_1^{-1} \nabla \cdot (A_1 \nabla u_2^{s}) - k^2 u_2^{s}, & \qquad \mbox{in} \,\,\, \Omega_{1} \nonumber \\
& \nonumber \\
J_2^{-1} \nabla \cdot (A_2 \nabla u_2^{s+1}) + k^2 u_2^{s+1} = \Psi_2(u_1^s), & \qquad \mbox{in} \,\,\, \Omega_{2} \label{eq:tran2} \\
\Psi_2(u_1^s) = - J_2^{-1} \nabla \cdot (A_2 \nabla u_1^{s}) - k^2 u_1^{s}, & \qquad \mbox{in} \,\,\, \Omega_{2} \nonumber 
\end{align}
where $\Psi_1$ and $\Psi_2$ are the source transfer function, $s$ is the iteration
step, $s = 0, 1, 2, \ldots$ 
Note that the transfered source $\Psi_1(u_2^s) = 0$ for $y < \tilde{l}$ or $y > \tilde{l} + \lpml$ ,
thus it has a compact support in the PML layer, so does $\Psi_1(u_2^s)$.
At last, the PML solutions on subdomains are summed up as the solution obtained by 
the domain decomposition method,
\begin{equation} \label{eq:ddm-sum}
  u_{\text{DDM}} = \sum_{s=0}^{\infty} (u_1^s + u_2^s).
\end{equation}
Although the PML equation \eqref{eq:dm1}-\eqref{eq:tran2} sovles the
truncated Helmholtz equation in the subdomain approximately, the
convergence of the series \eqref{eq:ddm-sum} to the solution of
\eqref{eq:pml} could be shown by

\begin{align*}
&\, \L \left(\sum_{s = 0}^N (u_1^s + u_2^s) \right) - f \\
       =&\, \L (u_1^0 + u_2^0) - f + \L \left(\sum_{s = 1}^N (u_1^s + u_2^s) \right) \nonumber \\
       =&\, -\Psi(u_1^0) - \Psi (u_2^0) + \L (u_1^1 + u_2^1) + \L \left(\sum_{s = 2}^N (u_1^s + u_2^s) \right) \nonumber \\
       =&\, -\Psi(u_1^1) - \Psi (u_2^1) + \L (u_1^2 + u_2^2) + \L \left(\sum_{s = 3}^N (u_1^s + u_2^s) \right) \nonumber \\
       =&\,  \ldots \nonumber \\
       =&\, \L ( u_1^N + u_2^N),  \nonumber
\end{align*}
and the remaining term $\L (u_1^N + u_2^N) \rightarrow 0$ as $N \rightarrow \infty$, which could be ensured
by the convergence of the PML method \cite{Chen2010} together with the analysis of wave traveling
in layered medium as follows.

%
%
%
%

\begin{figure}[ht!]
\centerline{
\includegraphics[width=.6\textwidth]{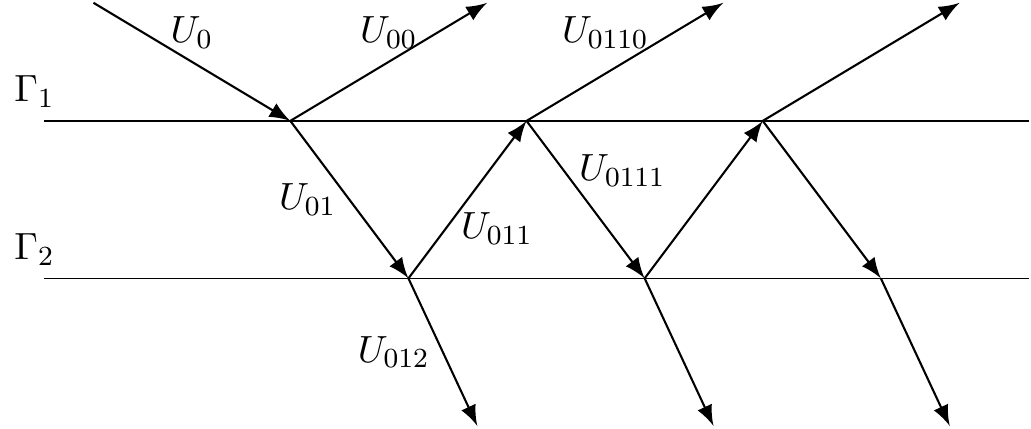}
} 
\caption{Wave traveling in three layered medium.} 
\label{fig:travel}
\end{figure}

The solution of the domain decomposition method 
in the form of \eqref{eq:ddm-sum} 
could be interpreted as the superposition 
of the incident waves, reflected waves and refracted waves 
that propagate in the layers \cite{Chew1999}, 
as is illustrated in Fig \ref{fig:travel}.

Suppose the incident wave $U_0$ comes from the upper layer, 
then at interface $\Gamma_1$, $U_0$ causes
a reflected wave $U_{00}$ going upwards in the upper layer 
and a refracted wave $U_{01}$ going downwards in the middle layer.
The wave $U_{0}+U_{00} + U_{01}$ is approximately the solution $u_1^0$ of the subdomain equation \eqref{eq:dm1} with $i = 1$.

Then at interface $\Gamma_2$, $U_{01}$ causes
a reflected wave $U_{012}$ going downwards in the lower layer 
and a refracted wave $U_{011}$ going upwards in the middle layer.
The wave $U_{01} + U_{012} + U_{011}$ is approximately the solution $u_2^1$ of the subdomain equation \eqref{eq:tran2}.

Then at interface $\Gamma_1$, $U_{011}$ causes
a reflected wave $U_{0110}$ going upwards in the upper layer 
and a refracted wave $U_{0111}$ going downwards in the middle layer.
The wave $U_{011} + U_{0110} + U_{0111}$ is approximately 
the solution $u_1^2$ of the subdomain equation \eqref{eq:tran1}.
The traveling process goes on, 
and the superposition of all the waves is the solution to \eqref{eq:pml},
\begin{align}  \label{eq:series}
u&= U_0 + U_{00} + U_{01} + U_{012} + U_{011} \nonumber \\
 & \quad + U_{0110} + U_{0111} + U_{01112} + U_{01111} + \dots
\end{align}
and the series \eqref{eq:series} is approximately the series \eqref{eq:ddm-sum}.


The convergence of the new overlapping domain decomposition method related 
closely to the medium perporty of the layers and the size 
of the overlapping region. 
When the overlapping region of the subdomains lies inside the middle
layer of the three, e.g., $\tilde{l} < d$, the convergence rate of the domain
decomposition method is at most the convergence rate of the series \eqref{eq:series}.
The worst case happens when there is a narrow wave guide, and 
the overlapping domain lies inside the wave guide,
e.g. $k_2 > k_1 = k_3$, $\tilde{l} < d$ and $d$ is small.
To avoid such cases, the overlapping region 
should have a non-zero minimum size.

 The overlapping region ensures the convergence of the new domain decomposition
  method for layered medium.
  The convergence of non-overlaping DDM might deteriorate if the subdomain
  interface lies right in a waveguide. We have two remarks on the new domain decomposition method. 

\textbf{Remark 1}: The convergence of the solution enables direct solving 
the Helmholtz
equation with the method, rather than use it as a preconditioner, 
which is crucial for our new fast method. 

\textbf{Remark 2}: An extend PML layer could be defined that it includes
a PML layer and a layer that doesn't absorb at all, for example, 
the layer $[-\lx, \lx] \times [0, \tilde{l}+\lpml]$ 
is an extend PML layer.
Since it's all about the PML layer parameters, we do not make a 
distinction between the two and simply call them the PML layer.

%
%
\subsection{Mapping instead of solving}

The domain decomposition method in the above subsection could be
revised that the solving of PML equation on subdomains 
\eqref{eq:tran1}-\eqref{eq:tran2} is substituted by mapping.

For subdomain $\Omega_1$, a mapping $\G_1$ from incidents trace $U^{\II}$ on the line $[-\lx, \lx]
\times 0$
to the wave solution $\bar{u}$  in $\Omega_1$ 
is defined as follows:
Given $U^{\II}$ on the line $[-\lx, \lx] \times 0$,
solve $\hat{u}$ as its extension such that 
\begin{align} \label{eq:Textn}
  J_2^{-1} \nabla \cdot (A_2 \nabla \hat{u}) + k^2 \hat{u} = 0, & \qquad \mbox{in} \,\,\, [-\lx,\lx] \times [0, \tilde{l} + \lpml]\\
  \hat{u} = U^\II, & \qquad \mbox{on} \,\,\, [-\lx,\lx] \times 0
\end{align}
It's obvious that if $U^\II$ is the trace of a solution to \eqref{eq:tran2},
then $\hat{u}$ is the restriction of that solution 
on the region $[-\lx,\lx] \times [0, \tilde{l} + \lpml]$.
The extension $\hat{u}$ is then transfered as source, 
\begin{equation} 
\Psi_1(\hat{u}) = - J_1^{-1} \nabla \cdot (A_1 \nabla \hat{u}) - k^2 \hat{u},  \qquad \mbox{in} \,\,\, \Omega_{1} \label{eq:Ttran}, 
\end{equation} 
with which the wave field solution $\bar{u}$ to PML equation in 
subomain $\Omega_1$ is solved
\begin{equation} 
J_1^{-1} \nabla \cdot (A_1 \nabla \bar{u}) + k^2 \bar{u}= \Psi_1(\hat{u}),  \qquad \mbox{in} \,\,\, \Omega_{1} \label{eq:Thelm} .
\end{equation} 
The mapping is then defined as $\bar{u} = \G_1(U^\II) $.

Another mapping $\F_1$ from incidents trace $U^{\II}$ on the line $[-\lx, \lx]
\times 0$ 
to the field trace $U^\FF$ on the same line, is defined by
$\DD U^\FF = \F_1(U^\II) \triangleq \left. \G_1(U^\II) \right|_{[-\lx, \lx] \times 0} $ . 
Although both the incident trace and the field trace is on the line $[-\lx, \lx] \times 0$ ,
it is referred as incident boundary or field boundary, respectfully. 
For subdomain $\Omega_2$, similar mapping $\G_2$ and $\F_2$ could be defined.

Now the domain decomposition method for Helmholtz equation with three
layered medium could be revised as follows:
first, solve the subdomain problem with the restricted source,   
\begin{equation} 
  J_i^{-1} \nabla \cdot (A_i \nabla u_i) + k^2 u_i = f_i, \qquad \mbox{in} \,\,\, \Omega_i, i = 1,2
\end{equation}
where $f_1 = f \cdot \chi_{y < 0}$ for $\Omega_1$, and $f_2 = f \cdot
\chi_{y \geq 0}$ for $\Omega_2$, the solution is denoted $u_i^0$ for
$i = 1,2$, and the field trace of the solutions are $U_i^{\FF,0} =
\left. u_i^{0} \right|_{[-\lx, \lx] \times 0}$, for $i = 1,2$.

Then each subdomain takes its neighbor's field trace as its own incident
trace, map the incident trace to filed trace, and so on,
\begin{align}
U_1^{\II,s+1} &= U_2^{\FF,s} &\qquad \mbox{in} \,\,\, \Omega_{1} \nonumber\\
U_1^{\FF,s+1} &= \F_1(U_1^{\II,s+1}) & \qquad \mbox{in} \,\,\, \Omega_{1} \nonumber\\
\\
U_2^{\II,s+1} &= U_1^{\FF,s} & \qquad \mbox{in} \,\,\, \Omega_{2} \nonumber\\
U_2^{\FF,s+1} &= \F_2(U_2^{\II,s+1}) &\qquad \mbox{in} \,\,\, \Omega_{2} \nonumber
\end{align}
for $s = 0, 1, 2, \ldots$, and the domain decomposition solution is 
\begin{equation} 
  u_{\text{DDM}} = u_1 + u_2 + \G_1 \left(\sum_{k=0}^{\infty} U_1^{\II,s} \right) 
   + \G_2 \left(\sum_{s=0}^{\infty} U_2^{\II,s} \right)   .
\end{equation}

%
%
\subsection{STDDM with four subdomains}

The above domain decomposition method with two subdomain in $y$ direction 
could  be easily
extended to four subdomains in both $x$ and $y$ directions. The major
difference is that the incident boundaries, field boundaries and their source
tranfer regions are a little complicated for four subdomains.

\begin{figure}[ht!] 
\includegraphics[width=\textwidth]{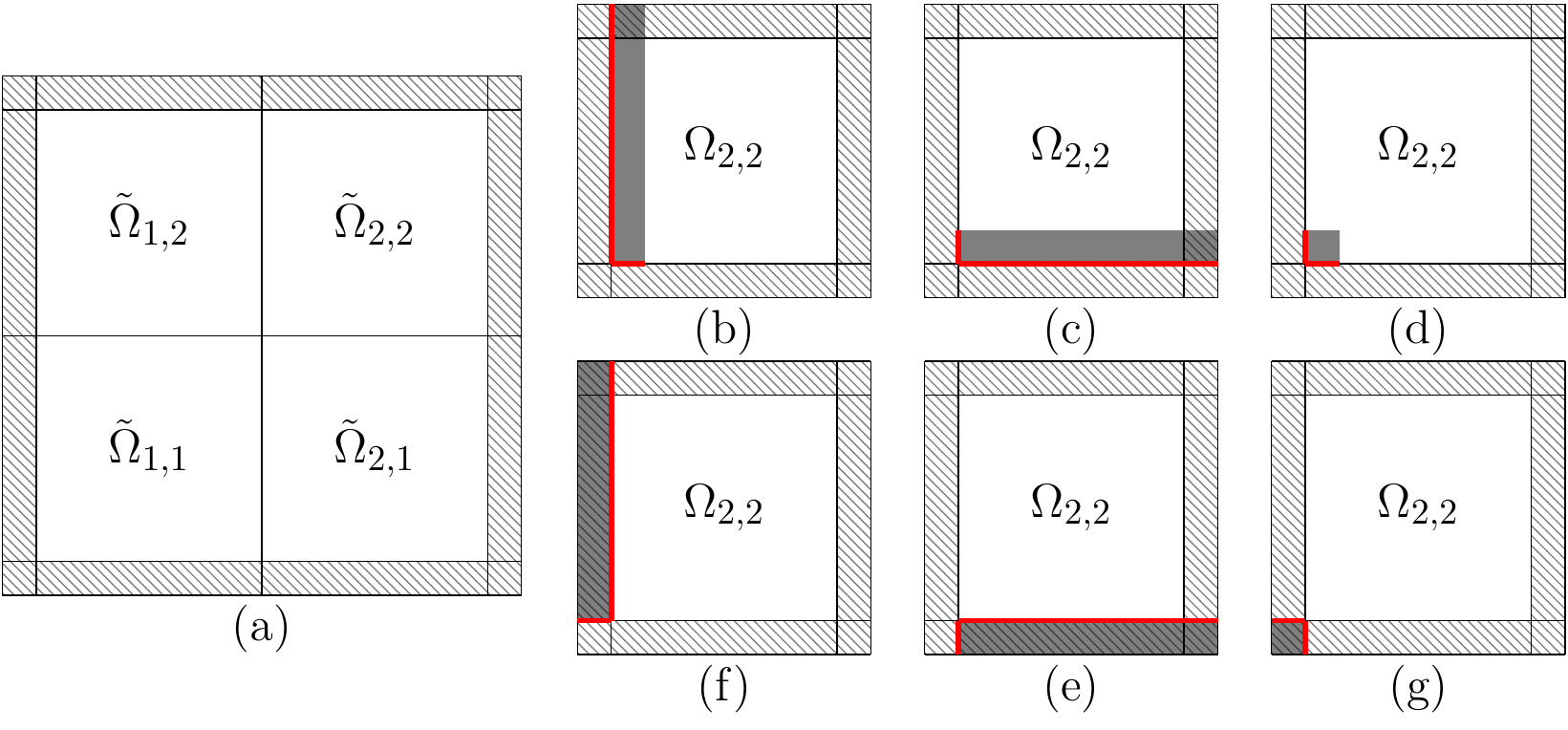} 
  
\caption{Domain decomposition with four subdomains. The hatched area 
  is the PML layer, the shaddowed area is the source
  transfer region, the thick lines are the incident or field
  boundaries. (a) four subdomain's interior region $\tilde{\Omega}_{i,j}$,
  $i,j=1,2$ and the PML layer of total domain. (b-d) incident boundaries and
  corresponding source transfer region of subdomain $\Omega_{2,2}$.
  (e-g) field boundaries and corresponding source transfer region of
  subdomain $\Omega_{2,2}$. }
\label{fig:fourdm}
\end{figure}

The total domain $\Omega$ is decomposed into four smaller subdomains
$\Omega_{i,j}$, $i,j=1,2$. The interior (region without PML layer) of the subdomain
$\Omega_{i,j}$ are denoted $\tilde{\Omega}_{i,j}$, they are non-overlapped 
and their union is the interior of the total domain,
as is shown in Fig \ref{fig:fourdm}\,-(a).  Each subdomain
$\Omega_{i,j}$ has its PML layer lie in its neighbors.

There are three kind of incident boundaries, denoted
$\Gamma_{i,j}^{\II}$, and three kind of field boundaries, denoted
$\Gamma_{i,j}^{\FF}$, for subdomain $\Omega_{i,j}$, as in Fig \ref{fig:fourdm}\,-(b-g). 
For examples, on subdomain $\Omega_{2,2}$, the incident boundary for
wave comes from subdomain $\Omega_{1,2}$ is shown in Fig
\ref{fig:fourdm}\,-(b), 
and the field boundary for wave goes to
subdomain $\Omega_{1,2}$ is shown in Fig \ref{fig:fourdm}\,-(e).

The  incident traces on boundary $\Gamma_{i,j}^\II$ are denoted as
$U_{i,j}^{\II}$, 
 and the field traces on boundary $\Gamma_{i,j}^\FF$ are denoted as
 $U_{i,j}^{\FF}$.
 The mapping from the incident trace to the solution on subdomain
 $\Omega_{i,j}$ is denoted $\G_{i,j}$,
 while the the mapping from the incident trace to the field trace on subdomain
 $\Omega_{i,j}$ is denoted $\F_{i,j}$

The domain decomposition method with four subdomains is shown in 
Algorithm \ref{alg:ddm4}.
In the algorithm, the wave propagates between children subdomains
via the iteration (3\,-7), we call it the iteration of incident and field traces 
 from now on.

%
%
%
%
%
%
\begin{breakablealgorithm}
\caption{Domain decomposition with four subdomains.} \label{alg:ddm4}
\begin{algorithmic}[1]

  \State Solve the mapping $\F_{i,j}$ on subdomain $\Omega_{i,j}$, $i, j= 1,2$, 
  \Statex[1] with direct solver.
 
  \State Solve the local problem on $\Omega_{i,j}$ with source 
  $f_{i,j} = f | \tilde{\Omega}_{i,j}$, 
  \Statex[1] restrict the solution
  $u_{i,j}^0$ to field trace $U_{i,j}^{\FF,0}$.

    \While{$\sum_{i,j=1,2}||U_{i,j}^{\FF,s}|| > \varepsilon$}
        \State Send subdomain $\Omega_{i,j}$'s field trace $U_{i,j}^{\FF,s}$
        \Statex to its siblings $\Omega_{i',j'}$ as incident trace $U_{i',j'}^{\II,s+1}$
        \State Record the incident traces $U_{i,j}^{\II,s+1}$

        \State Map the incidents trace to field trace
               $U_{i,j}^{\FF,s+1} = \F_{i,j}\FSUP (U_{i,j}^{\II,s+1}) $
        \State Set $s = s+1$
    \EndWhile

    \State Solve the local problem 
    on $\Omega_{i,j}$ with the summation
    of incident traces using direct solver, 
     the solution is denoted  $\G_{i,j} \left(\sum_{s > 0} U_{i,j}^{\II,s}\right)$.

    \State Sum up the solutions of all subdomains to get the total solution
    $$u = \sum_{i,j=1,2} \left(u_{i,j}^0 + \G_{i,j} \left(\sum_{s > 0} U_{i,j}^{\II,s}\right) \right). $$

\end{algorithmic}
\label{alg-src-up} 
\end{breakablealgorithm}

%
%
%
%
\section{The Fast Propagation Method}

%
%
\subsection{Hierarchical domain decomposition}

\begin{figure}[ht!] \label{Fig-hddm}
\begin{center}

\includegraphics[width=.22\textwidth]{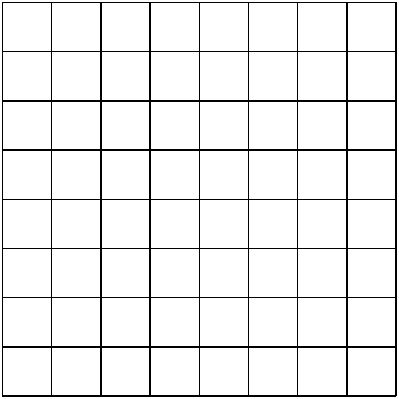} 
\includegraphics[width=.22\textwidth]{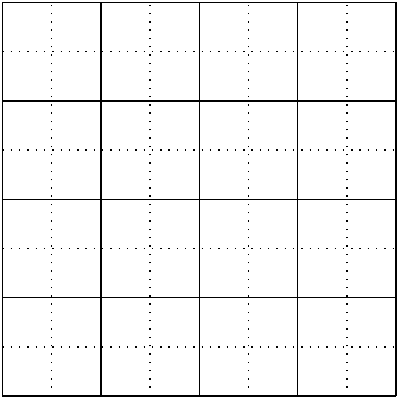} 
\includegraphics[width=.22\textwidth]{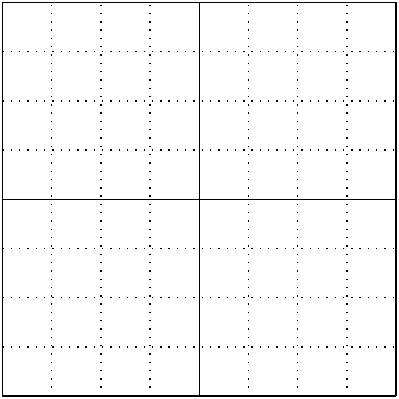} 
\includegraphics[width=.22\textwidth]{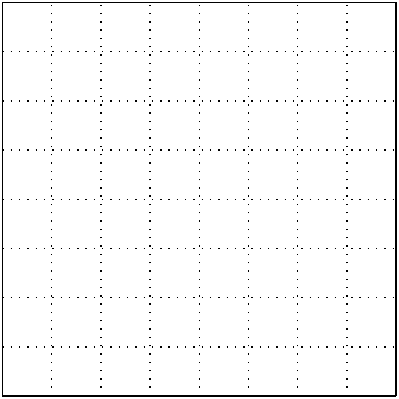} 

\end{center} 
\caption{Hierarchical Domain decomposition with 4 levels.
  From left to right, domain decomposition on level 1, 2, 3, 4.}
\end{figure}

A rectangular domain of $[-L,L]$ is decomposed into smaller
rectangular blocks (or subdomain) on different levels.
Denote the number of levels as $N_L + 1$, 
the level $l = 0$ is referred as the bottom level 
and the level $l = N_L$ is referred as the top level.
The number of blocks in $x$
direction at level $l$ is $2^{N_L-l}$, where $l = 0, \ldots, N_L$.
Let $I_l = \left\{1, \ldots, 2^{N_L-l} \right\}$,
on the level $l$, the block which is the $i$-th block in $x$ direction
and the $j$-th block in $y$ direction, is denoted $\Omega_{i,j;l}$,
 where $i,j \in I_l$.
Each block shares an overlapping PML layer region of length $\lpml$
with its neighbors on the same level.

The quadtree structure of the multiple level domain decomposition is 
built as follows.
Each block $\Omega_{i,j;l}$ on level $l = 2^{N_L-l}, \dots, 1$ has four children
$\DD \Omega_{2i-1,2j-1;\l-1}$, 
$\DD \Omega_{2i-1,2j;\l-1}$, 
$\DD \Omega_{2i,2j-1;\l-1}$, 
and $\DD \Omega_{2i,2j;\l-1}$
on level $l-1$.  
For simplicity, the children of block $\Omega_{i,j;\l}$
is denoted $\Omega_{i',j';\l-1}$, where $i' = 2i-1, 2i$, $j' = 2j-1, 2j$.
On the other hand, each block $\DD \Omega_{i,j;\l}$ on level $l$
has a father $\DD \Omega_{\lceil i/2\rceil, \lceil j/2
  \rceil;\l+1}$ on level $l + 1$, where $l < N_L$.
The father-son relationship of the blocks leads to the quadtree structure.

\begin{figure}[ht!] 
\begin{center}
\includegraphics[width=\textwidth]{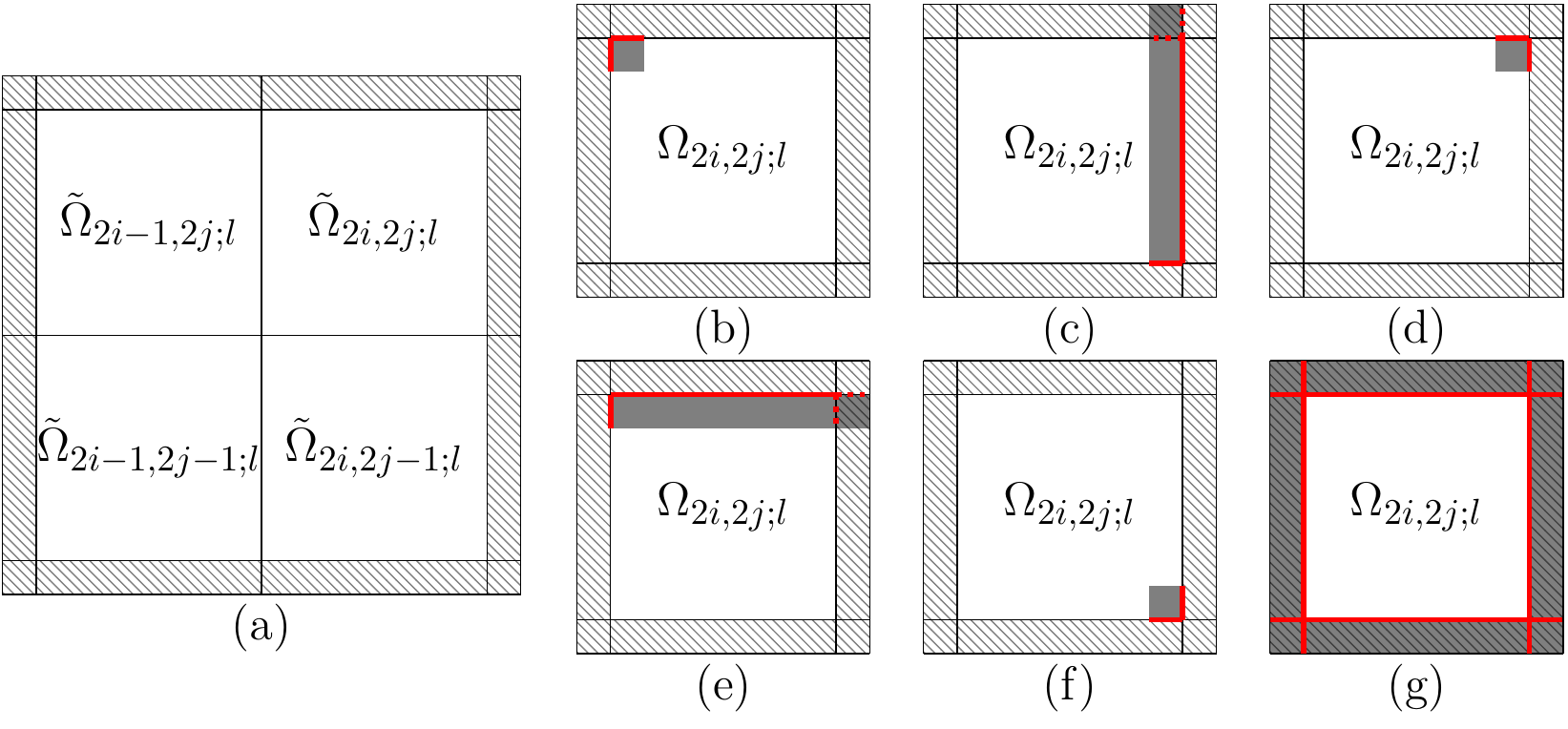} 

\end{center} 
\caption{Incident boundaries and field boundaries extension. 
  The hatched area
  is the PML layer, the shaddowed area is the
  source transfer region, the thick lines are the incident or field
  boundaries. (a) four children block's interior region and the PML layer of
  father block $\Omega_{i,j;\l+1}$.  
  (b-f) ascendant's incident boundaries and corresponding
  source transfer region on child $\Omega_{2i,2j;\l}$. (e) total field
  traces and corresponding source transfer region on child $\Omega_{2i,2j;\l}$. }
\label{fig:fourdm2}
\end{figure}

The incident boundaries and field boundaries on block $\Omega_{i,j;\l}$ 
include not only the boundaries between siblings as in Fig \ref{fig:fourdm}, 
but also its ascendant's incident boundaries and field boundaries, 
as is shown in Fig \ref{fig:fourdm2}. 
We call the boundaries as in Fig \ref{fig:fourdm}
 the corresponding incident and field boundaries 
between siblings.
The incident boundary of block $\Omega_{i,j;\l}$ is denoted
$\Gamma_{i,j;\l}^{\II}$, and the field boundary of block
$\Omega_{i,j;\l}$ is denoted $\Gamma_{i,j;\l}^{\FF}$.
We see 
$\DD \Gamma_{i,j;l}^{\II}  \subset \mathop{\cup}_{i',j'} \Gamma_{i',j';l-1}^{\II} $
and 
$\DD \Gamma_{i,j;l}^{\FF}  \subset \mathop{\cup}_{i',j'} \Gamma_{i',j';l-1}^{\FF} $.
The mapping from incident trace to solution on the block 
is denoted $\G_{i,j;l}$, 
and the mapping from incident trace to field trace 
on the block is denoted $\F_{i,j;l}$ .

%
%
\subsection{Setup phase}

In the setup phase, the mapping from incident traces to field traces
is constructed bottom up level by level.
The mapping on block $\Omega_{i,j}^{\,0}$, $i,j \in I_0$ could be computed
with external direct solver, and the mapping on block $\Omega_{i,j;\l}$
of level $l > 0$ is computed as follows. 

Given an incident $\delta$ lies in $\Gamma_{i,j;\l}^{\II}$ , it must 
lie in the incident boundary of one of the children, denoted as
$\Omega_{i_0,j_0;\l-1}$.  First the local problem on children
$\Omega_{i_0,j_0;\l-1}$ with source $\delta$ is considered, and the
field trace of the solution on $\Gamma_{i_0,j_0;\l-1}^{\FF}$ is solved
by mapping $\F_{i_0,j_0;\l-1}\FSUP$.  Then the field trace of
$\Omega_{i_0,j_0;\l-1}$ is send to its siblings as incidents, 
and the iteration of incident and filed trace between siblings applies,
and the incident trace in the iteration is denoted $U_{i',j';\l-1}^{\II,s}$, 
where $s$ is the iteration number. 
At last,
field trace on $\Gamma_{i',j';\l-1}^{\FF}$ caused by sum of incidents
computed with the mapping $\F_{i',j';\l-1}$, 
along with the field trace caused by $\delta$ on
$\Omega_{i_0,j_0;\l-1}$, are add up as the field trace 
$U_{i,j;\l}^{\FF}$ on $\Omega_{i,j;\l}$
 caused by $\delta$,
\begin{equation}
U_{i,j;\l}^{\FF} = 
     \F_{i_0,j_0;\l-1}\FSUP  (\delta) \Big|_{\Gamma_{i,j;\l}^{\FF}}
    + \sum_{i',j'} \F_{i',j';\l-1}\FSUP (\sum_{s} U_{i',j';\l-1}^{\II,s}) 
    \Big|_{\Gamma_{i,j;\l}^{\FF}},
\end{equation}
and the mapping is 
\begin{equation}
\F_{i,j;\l}(\delta) = U_{i,j;\l}^{\FF}.
\end{equation}

The algorithm of building the mapping from incident traces to 
field traces is as follows.
\begin{breakablealgorithm}
\caption{Build mapping of incident traces to field traces} \label{alg:mapp} 
\begin{algorithmic}[1]

\State On level 0, build the mapping of incident to field with direct solver.

\For {levels $l = 1, \ldots, N_L$}

    \State On block $\Omega_{i,j;\l}$, 
    \For {incident $\delta $ lies in $ \Gamma_{i,j;\l}^{\II}$ }

       \State Find the children $\Omega_{i_0,j_0;\l-1}$ such that
              $\delta$ lies in $ \Gamma_{i_0,j_0;\l-1}^{\II}$
       \State On children $\Omega_{i_0,j_0;\l-1}$,
       \Statex map the incidents $\delta$ to
               field trace $U_{i_0,j_0;\l-1}^{\FF}$, 
       \Statex and add part of them to father's
               field trace  $U_{i,j;\l}^{\FF}$.

       \State Set $U_{i_0,j_0;\l-1}^{\FF,0} =
              U_{i_0,j_0;\l-1}^{\FF}$ on children
              $\Omega_{i_0,j_0;\l-1}$,
       \Statex and set $U_{i',j';\l-1}^{\FF,0} = 0$ on other children $\Omega_{i',j';\l-1}$.
       \While{$\sum_{i',j'} ||U_{i',j';\l-1}^{\FF,s}|| > \varepsilon$}

           \State Send the children's corresponding field trace
                  $U_{i',j';\l-1}^{\FF,s}$ 
           \Statex[4] to its sibilings as incidens
                  $U_{i',j';\l-1}^{\II,s+1}$

           \State Map the incidents to field trace
                  $U_{i',j';\l-1}^{\FF,s+1} = \F_{i',j';\l-1}\FSUP (U_{i',j';\l-1}^{\II,s+1}) $

           \State Set $s = s+1$
       \EndWhile
       \State Map the sum of incidents to field trace on children, 
       \Statex and add them to father's field $U_{i,j;\l}^{\FF}$.

    \EndFor

\EndFor

\end{algorithmic}
\end{breakablealgorithm}

%
%
\subsection{Solve phase}

With the mapping of incident traces to filed traces that is 
constructed on each block
of all levels, the Helmholtz equation could be solved in two phases,
the source-up phase and the the solution-down phase.


\subsubsection{The Source-up phase}
In the source-up phase, the wave propagates on all levels bottom up 
as incident traces.

The following problem is considered,
for the block $\Omega_{i,j;\l}$, the local
solution on its four children are known, 
e.g., $u_{i',j';\l-1}^{0}$,
so does their field traces 
$U_{i',j';\l-1}^{\FF,0}$,
how to solve
the solution $u_{i,j;\l}$ on $\Omega_{i,j;\l}$, 
and its field trace $U_{i,j;\l}^{\FF}$. 
The iteration of incident and filed trace between siblings applies directly,
denote the incident traces in the iteration as $U_{i',j';\l-1}^{\II, s}$,
and
the solution on $\Omega_{i,j;\l}$ is
\begin{equation}
\DD u_{i,j;\l} = \sum_{i',j'} \left( u_{i',j';\l-1}^{0} 
          + \G_{i',j';\l-1}(\sum_{s} U_{i',j';\l-1}^{\II,s}) \right),
\end{equation}
and the the field trace of $u_{i,j}^{l}$ is
\begin{align}
\DD U_{i,j;\l}^{\FF} &=  \sum_{i',j'} \left( 
          U_{i',j';\l-1}^{\FF,0} \Big|_{\Gamma_{i,j;\l}^{\FF}} 
          + \F_{i',j';\l-1}\FSUP(\sum_{s} U_{i',j';\l-1}^{\II,s}) \Big|_{\Gamma_{i,j;\l}^{\FF}}
         \right) .
\end{align}

Review the procedure we found that to apply the procedure to next level,
the incident to field mapping operation $\F_{i',j';\l-1}$ of children is needed, 
while the solving operation 
$\G_{i',j';\l-1}(\sum_{s} U_{i',j';\l-1}^{\II,s})$  could be post processed.
Apply the procedure from bottom level to top level leads to the following
source-up algorithm.

\begin{breakablealgorithm}
\caption{Source-up} \label{alg:src-up} 
\begin{algorithmic}[1]

\Require {Right hand side $f$ of the linear system}
\Ensure {Solution $u_{i,j,0}^{0}$ on $\Omega_{i,j}^{0}$,
  \Statex[1] and sum of incidents on $\Omega_{i,j;\l}$ on level $l>0$}

\State On level $l=0$,
\Statex[1]  solve the local problem on $\Omega_{i,j}^{0}$ with the
  source $f_{i,j}^{0} = f | \hat{\Omega}_{i,j}^{0}$,
\Statex[1]  the solution  $u_{i,j,0}^{0}$ and the its field trace $U_{i,j,0}^{\FF}$ are recored.

\For {levels $l = 1, \ldots, N_L$}

    \State On block $\Omega_{i,j;\l}$, 
    \Statex[2] use the field trace $U_{i',j';\l-1}^{\FF}$ of
    the four childrens $\Omega_{i',j';\l-1}$,
    \Statex[2] add part of $U_{i',j';\l-1}^{\FF}$ to $U_{i,j;\l}^{\FF}$,     
    \Statex[2] set $U_{i',j';\l-1}^{\FF,0} = U_{i',j';\l-1}^{\FF}$, 
    \While{$||U_{i',j';\l-1}^{\FF,s}|| > \varepsilon$}

        \State Send children's corresponding field traces $U_{i',j';\l-1}^{\FF,s}$
        \Statex to its sibilings as incidens $U_{i',j';\l-1}^{\II,s+1}$

        \State Map the incidents to field 
               $U_{i',j';\l-1}^{\FF,s+1} = \F_{i',j';\l-1}\FSUP (U_{i',j';\l-1}^{\II,s+1}) $

        \State Set $s = s+1$
    \EndWhile

    \State Sum up the incidents $\sum_{s} U_{i',j';\l-1}^{\II,s}$ for childrens

    \Statex[2] and map them to the field trace $\Omega_{i,j;\l}$,
    \Statex[2] then add to  $U_{i,j;\l}^{\FF}$.
\EndFor

\end{algorithmic}

\end{breakablealgorithm}

The solution to the total problem could then be expressed as

\begin{equation} \label{eq:src-up}
\DD u = \sum_{i,j \in I_{0}} u_{i,j;0}^{0}
      + \sum_{l > 0}\sum_{i,j \in I_{l}} \sum_{i',j'}
      \left( \G_{i',j';\l-1}(\sum_{s} U_{i',j';\l-1}^{\II,s}) \right).
\end{equation}

\subsubsection{The Solution-down phase}
In the solution-down phase the wave propagates on all levels top down as 
incident traces. 

The solution \eqref{eq:src-up} resulting from Algorithm \ref{alg:src-up} 
still needs to solve the local Helmholtz problem  
with given incidents on blocks of different levels,
fortunately, the local solutions could be break down to lower and lower level
 till level 0.
We consider the following problem:
on the block
$\Omega_{i,j;\l}$, given the incidents
$\U_{i,j;\l}^{\II}$, how to solve 
$\G_{i,j;\l}(\U_{i,j;\l}^{\II})$.

First the incident traces
$\U_{i,j;\l}^{\II}$ is divided into the incident traces on children
$\U_{i,j;\l}^{\II} = \DD \sum_{i',j'}
\U_{i,j;\l-1}^{\II,0}$, then with the incident to field mapping
$\F_{i',j';\l-1}\FSUP$ on each children, field trace of children is
generated, e.g., $\U_{i',j';\l-1}^{\FF,0}$, 
then the iteration of incident and filed trace between siblings applies,
and the incident traces in the iteration is denoted as 
$\tilde{U}_{i',j';\l-1}^{\II, s}$.
At last the solution on $\Omega_{i,j;\l}$ is
\begin{equation}
\G_{i,j;\l} (\U_{i,j;\l}^\II)
= \sum_{i',j'} \left( \G_{i',j';\l-1} (\sum_{s} \U^{\II,s}_{i',j';\l-1}) \right)
\end{equation}

Apply the procedure from level $l=N_L$ to $l=1$, since there are already
sum of incidents $\sum_{s} U^{\II,s}_{i',j';\l-1}$ on children blocks
$\Omega_{i',j';\l-1}$, the incidents $\sum_{s} \U^{\II,s}_{i',j';\l-1}$
from $\Omega_{i,j;\l}$ should be added on children.
The algoritm is discribed as follows.

\begin{breakablealgorithm}
\caption{Solution down} \label{alg:sol-dn} 
\begin{algorithmic}[1]

\Require {Solution $u_{i,j;1}^{0}$ on $\Omega_{i,j;1}$,
  \Statex[1]  and sum of incidents on $\Omega_{i,j;\l}$ on level $l>1$}
\Ensure {Solution $u$ of the linear system}

\For {levels $l = N_L, \ldots, 1$}

    \State On block $\Omega_{i,j;\l}$, divide the sum of incidents to its children,
    $$ \tilde{U}_{i,j;\l}^{\II} = \DD \sum_{i',j'} \U_{i',j';\l-1}^{\II,0}$$

    \State Map the incidents to field 
           $\U_{i',j';\l-1}^{\FF,1} 
           = \F_{i',j';\l-1} (\U_{i',j';\l-1}^{\II,0}) $

    \While{$||\tilde{U}_{i',j';\l-1}^{\II,s}|| > \varepsilon$}

        \State Send children's corresponding field trace $\U_{i',j';\l-1}^{\FF,s}$
        \Statex to its siblings as incidents $\U_{i',j';\l-1}^{\II,s+1}$

        \State Map the incidents to field 
               $\U_{i',j';\l-1}^{\FF,s+1} = \F_{i',j';\l-1}\FSUP (\U_{i',j';\l-1}^{\II,s+1}) $


        \State Set $s = s+1$
    \EndWhile

    \State Add the sum of incidents on children $\Omega_{i',j';\l-1}$,
       $$ \tilde{U}_{i',j';\l-1}^{\II} := \sum_{s} U_{i',j';\l-1}^{\II,s} + \sum_{s} \U_{i',j';\l-1}^{\II,s}$$ 

\EndFor

\State On level $l=0$,
\Statex[1]  solve the local problem on $\Omega_{i,j}^{0}$ with the incidents $\tilde{U}_{i,j;0}^{\II}$,
\Statex[1]  and add the solution to total solution $u$.

\end{algorithmic}
\end{breakablealgorithm}

Now the solution to the total problem is
\begin{equation} \label{eq:sol-dn}
\DD u = \sum_{i,j \in I_{0}} 
     \left(  u_{i,j;0}^{0}
      + 
      \G_{i,j;0}(\sum_{s} \tilde{U}_{i,j;\,0}^{\II,s})
      \right).
\end{equation}

%
%
%
%
\section{Numerical experiments}

\begin{figure}[ht!] 
\begin{center}
\includegraphics[width=.9\textwidth]{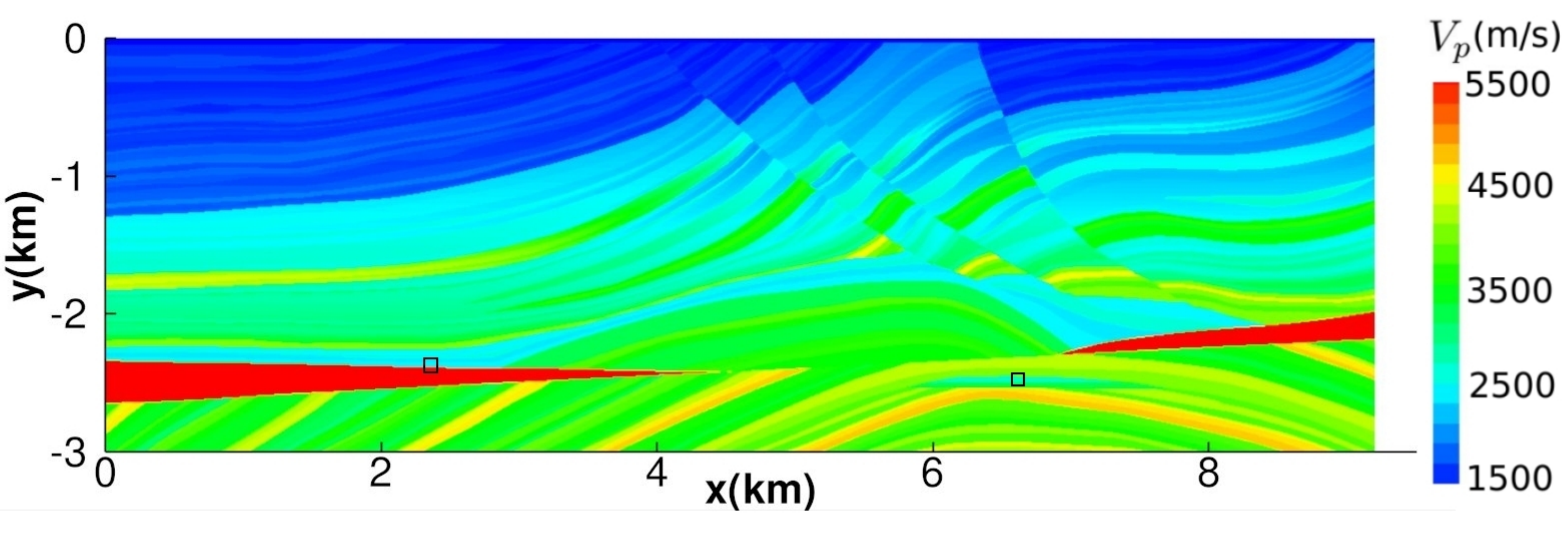} 
\end{center}
\caption{Velocity profile of Marmousi model. 
The solution with $N_{L} = 5$ is sampled in two small boxes in the figure.}
\label{fig:mar-vel}
\end{figure}

The new method is tested on the 2D Marmousi model in seismology,
which is  $3,000$ m deep and $9,200$ m wide.  Only P-wave
is considered, thus elastic wave equation becomes an acoustic equation.
The velocity profile is shown in Fig \ref{fig:mar-vel}, the
maximum velocity is 5500 km/s and the minmum velocity is 1500 km/s. 

\begin{figure}[ht!] 
\begin{center}
\includegraphics[width=.4\textwidth]{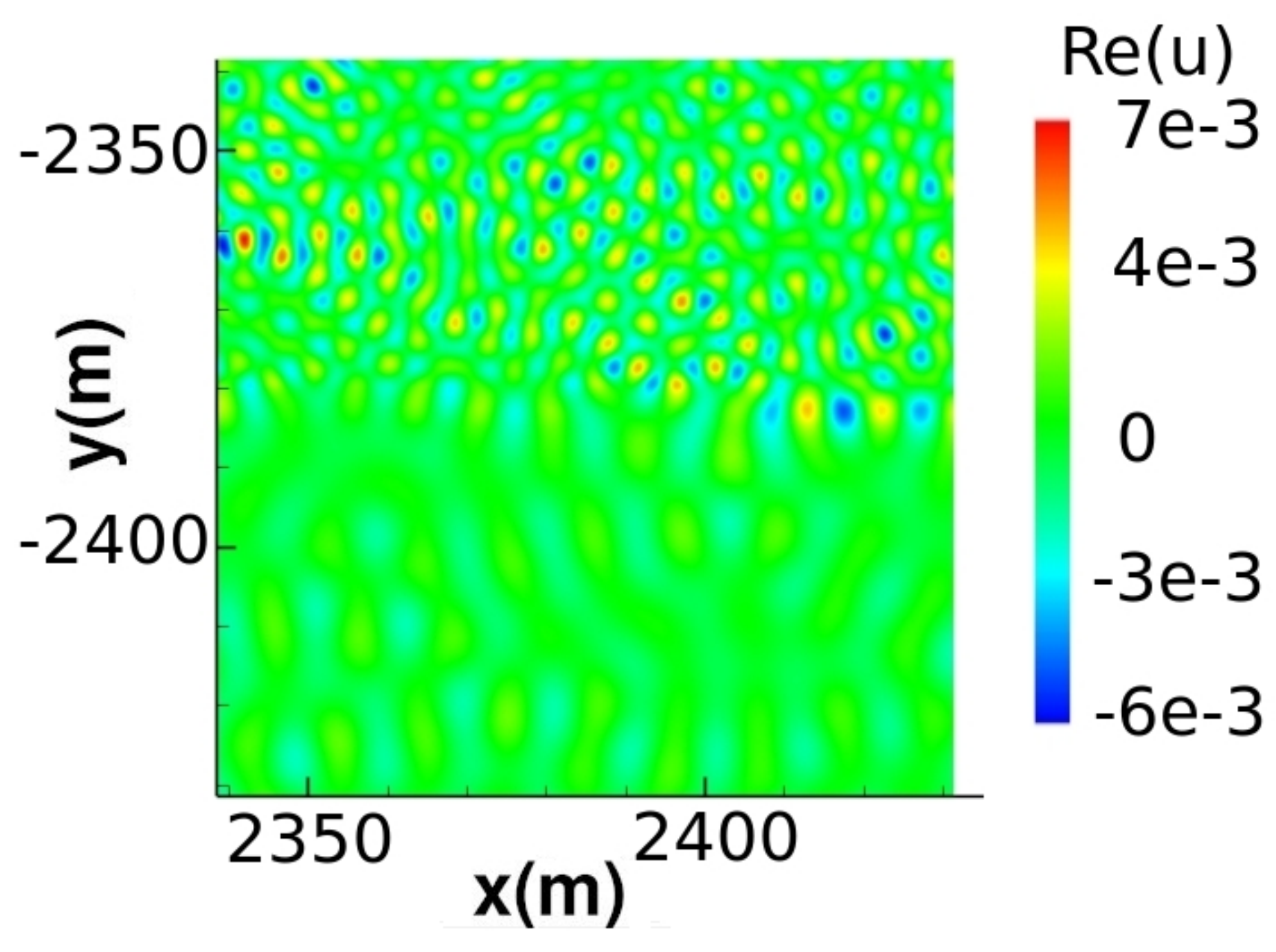} \qquad
\includegraphics[width=.4\textwidth]{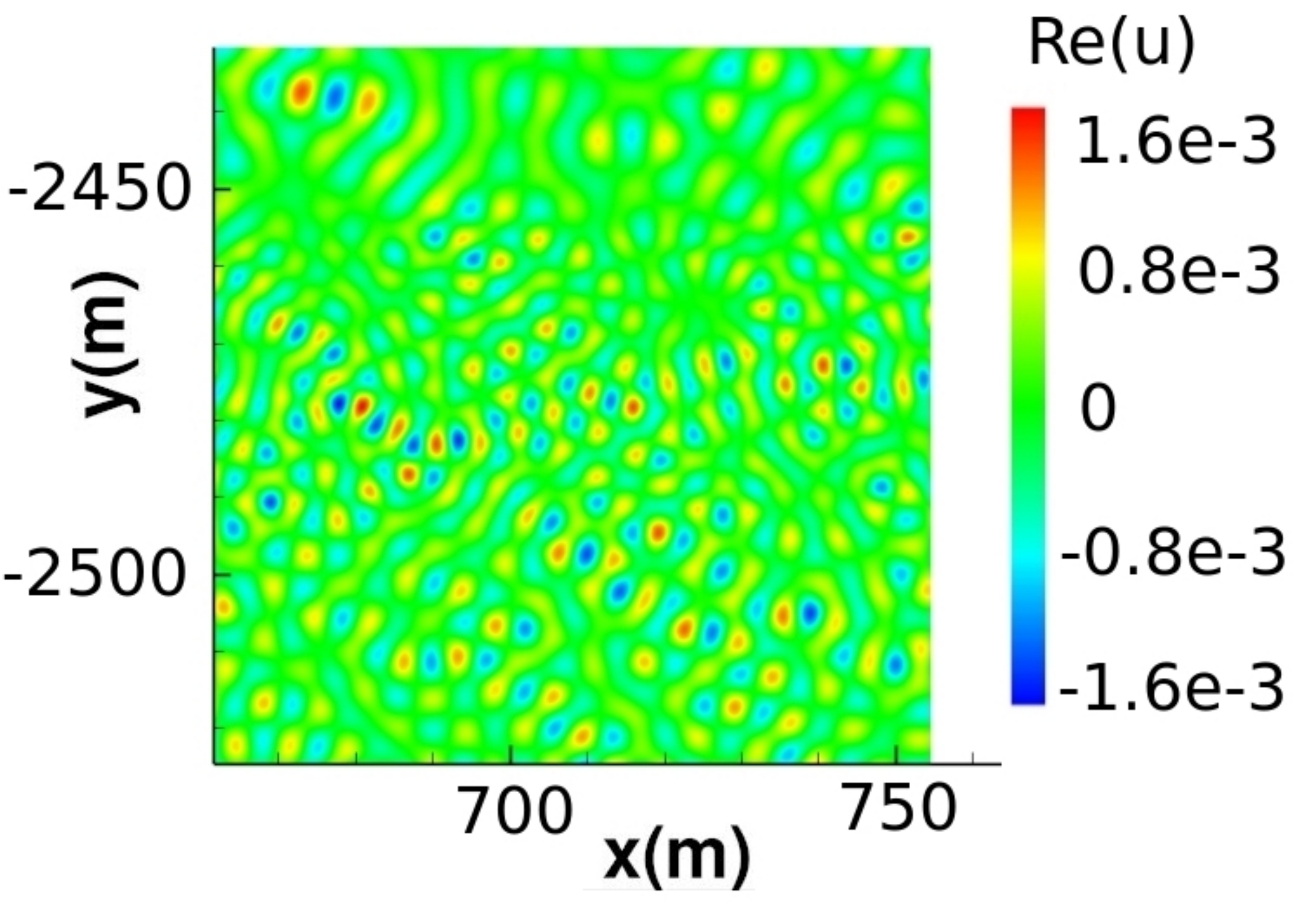} 
\end{center}
\caption{Real part of the solution with $N_{L} = 5$ in two small boxes
as marked in Fig \ref{fig:mar-vel}.}
\label{fig:mar-sol}
\end{figure}

Finite difference method with second order of accuracy is used to
discretize the Helmholtz equation. The block size on bottom level is 400
$\times$ 400,
ant the PML layer is of 40 grid points width.
Single shot in the corner of the domain 
at $(400 h_x, 400 h_y)$ is taken as the
source, where $h_x$, $h_y$
are the grid size in $x$ and $y$ direction, respectfully.
 The shape of the shot is an approximate delta function, 
 $\DD f_{i,j} = \frac{1}{h_x h_y}\delta(i - 400 h_x, j - 400 h_y)$.

\begin{table}[htbp]
\centering 
\begin{tabular}{|r|c|r|r|r|r|r|}
  \hline
  $N_L$   & Size                   & Freq            & No.   & Time  & Time  & Time  \\
          &                        & $\omega / 2\pi$ & procs & setup & solve & total \\
  \hline
  1     & 2,400 $\times$ 800     & 37              & 12    & 40    & 195   & 235   \\
  2     & 4,800 $\times$ 1,600   & 70              & 48    & 140   & 205   & 345   \\
  3     & 9,600 $\times$ 3,200   & 137             & 192   & 333   & 309   & 642   \\
  4     & 19,200 $\times$ 6,400  & 270             & 768   & 1212  & 685   & 1897  \\
  5     & 38,400 $\times$ 12,800 & 537             & 3,072 & 2891  & 883   & 3774  \\
  \hline
\end{tabular}
\caption{Time cost (in seconds) of the new method. } \label{tab:time}
\end{table}

\begin{table}[htbp]
\centering 
\begin{tabular}{|r|c|c|c|c|c|c|}
\hline
$N_L$ & Time    & Time    & Time    & Time    & Time    \\
      & Level 0 & Level 1 & Level 2 & Level 3 & Level 4 \\
\hline
1     & 39.6    & -       & -       & -       & -       \\
2     & 100     & 40.1    & -       & -       & -       \\
3     & 119     & 86.2    & 128     & -       & -       \\
4     & 127     & 74.7    & 256     & 754     & -       \\
5     & 129     & 98.7    & 355     & 716     & 1592    \\
\hline
\end{tabular}
\caption{Detailed setup phase time cost (in seconds).} \label{tab:setup}
\end{table}

The fast propagation method is suitable for parallel computing,
and could be easily extend to thousands of cores.
We test the method with different grid levels and grid sizes 
on cluster,
as listed in Table \ref{tab:time}. The tolerance of residual 
$\DD \frac{||A x - b||_2} {|| b ||_2} $ is $10^{-7}$.
Fig \ref{fig:mar-sol} shows the solution with $N_L = 5$ 
in two small boxes of $400 \times 400$ grid points 
as marked in Fig \ref{fig:mar-vel}.

The time cost of solving Helmholtz equation with the fast method in parallel  
is shown in Table \ref{tab:time}. The setup phase is the most demanding part 
in solving, since its complexity is $O(N^{3/2} \log N)$. The detailed time
cost in setup phase is shown in Table \ref{tab:setup}. 
The mapping on the bottom block is solved with direct solver, 
e.g. MUMPS \cite{MUMPS}, and 
the time cost is almost constant, since the bottom 
level block is of fixed size.
However, the time cost of building mapping on level $l+1$ 
is roughly twice of level $l$, where $l > 0$,  
which is time consuming for large Helmholtz problems.

%
%
%
%
\section{Conclusions}

A fast method is proposed for solving Helmholtz equations, the new method 
has a setup phase of complexity $O(N^{3/2} \log N)$ and a solve phase 
of complexity $O(N \log N)$. Our future work is to reduce the computation 
time of the new method by exploiting the low rank structure of the mappings 
and accelerating dense matrix operations with GPU.

\section*{Acknowledgments}
This work is supported by the National 863 Project of China
under the grant number 2012AA01A309, and the National Center for
Mathematics and Interdisciplinary Sciences of the Chinese Academy of
Sciences.


\end{document}